\documentclass{amsart}

\usepackage{amsthm,amssymb,amsmath}
\address{Steklov Mathematical Institute, Moscow\newline${}\quad$ Moscow State University
\newline${}\quad$ Kharkevich Institute for Information Transmission Problems, Moscow}

\email{agaif@mi.ras.ru}

\newcommand{\Z}{\mathbb{Z}}

\newcommand{\R}{\mathbb{R}}
\newcommand{\F}{\mathbf{F}}

\newcommand{\CC}{\mathcal{C}}
\newcommand{\RZ}{\mathcal{R}}
\newcommand{\hM}{\widehat{M}}
\newcommand{\Lob}{\mathbb{H}}

\newtheorem{theorem}{Theorem}[section]
\newtheorem{corollary}[theorem]{Corollary}

\newtheorem{proposition}[theorem]{Proposition}
\newtheorem{conjecture}[theorem]{Conjecture}
\newtheorem{question}[theorem]{Question}

\newtheorem{problem}[theorem]{Problem}

\newtheorem*{URC}{Universal Realisation of Cycles Property}

\theoremstyle{definition}

\newcommand{\link}{\operatorname{link}}
\newcommand{\codim}{\operatorname{codim}}
\newcommand{\Aut}{\operatorname{Aut}}
\newcommand{\Isom}{\operatorname{Isom}}

\title[Combinatorial realisation of cycles and small covers]{Combinatorial realisation of cycles and small covers}

\author[Alexander A. Gaifullin]{Alexander A. Gaifullin}

\thanks{The work was partially supported by RFBR  (projects 11-01-00694, 12-01-92104), by a grant of the President of Russian Federation (project MD-4458.2012.1), by a grant of the Government of the Russian Federation (project 2010-220-01-077), and by Dmitri Zimin's ``Dynasty'' foundation.}

\begin{document}

\begin{abstract}
In 1940s Steenrod asked if every  homology class $z\in H_n(X,\mathbb{Z})$ of every topological space~$X$ can be realised by an image of the fundamental class of an oriented closed smooth manifold. Thom found a non-realisable $7$-dimensional class and proved that for every~$n$, there is a positive integer~$k(n)$ such that the class $k(n)z$ is always realisable. The proof was by methods of algebraic topology and gave no information on the topology the manifold which realises the homology class. We give a purely combinatorial construction of a manifold that realises a multiple of a given homology class. For every~$n$, this construction yields a manifold~$M^n_0$ with the following universality property: \textit{For any~$X$ and $z\in H_n(X,\mathbb{Z})$, a multiple of~$z$ can be realised by an image of a (non-ramified) finite-sheeted covering of~$M^n_0$.\/} Manifolds satisfying this property are called URC-manifolds.
The manifold~$M^n_0$ is a so-called small cover of the permutahedron, i.e., a manifold glued in a special way out of $2^n$ permutahedra. (The permutahedron is a special convex polytope with $(n+1)!$ vertices.) Among small covers over other simple polytopes, we find a broad class of examples of URC-manifolds. In particular, in dimension~$4$, we find a hyperbolic URC-manifold. Thus we obtain that a multiple of every homology class can be realised by an image of a hyperbolic manifold, which was conjectured by Kotschick and L\"oh. Finally, we investigate the relationship between URC-manifolds and simplicial volume.

\end{abstract}

\maketitle

\section{Steenrod's problem on realisation of cycles}

\begin{problem}[Steenrod, 1940s]
For a given homology class $z\in H_n(X,\Z)$, does there exist an oriented smooth manifold~$M^n$ and a continuous mapping $f\colon M^n\to X$ such that $f_*[M^n]=z$? 
\end{problem}

This is known as Steenrod's problem on realisation of cycles.
The same question was asked for modulo~$2$ homology classes without the assumption of the orientability of~$M^n$. A homology class that can be realised by an image of the fundamental class of a smooth manifold will be called \textit{realisable\/}. We shall say that an integral homology class~$z$ is \textit{realisable with multiplicity~$k$\/} if the class~$kz$ is realisable.

\begin{theorem}[Thom~\cite{Tho54}, 1954]\label{theorem_Thom}
1) Every modulo~$2$ homology class is realisable.

2) For $n\le 6$, every $n$-dimensional integral homology class is realisable.

3) There exists a non-realisable $7$-dimensional integral homology class.

4) For each~$n$, there exists a positive integer $k=k(n)$ such that for every homology class $z\in H_n(X,\Z)$ of every topological space~$X$, the class~$kz$ is realisable.
\end{theorem}

Next important result was obtained by Novikov~\cite{Nov62} using the famous Milnor--Novikov computation of the rings of unitary and oriented cobordisms. He proved that the obstructions to the realisability of an $n$-dimensional homology class have odd order, and for every odd prime~$p$, the obstructions of order~$p$ lie in homology groups of dimensions $n-2q(p-1)-1$ for $q\ge 1$. In particular, he proved 

\begin{theorem}[Novikov~\cite{Nov62}, 1962]\label{theorem_Novikov}
Suppose that the group $H_{n-2q(p-1)-1}(X,\Z)$ has no $p$-torsion for every odd prime~$p$ and every~$q\ge 1$. Then any homology class $z\in H_n(X,\Z)$ is realisable.
\end{theorem}

Buchstaber~\cite{Buc69} related the obstructions to the realisability of homology classes to the differentials of the Atiyah--Hirzebruch spectral sequence. Recall that the Atiyah--Hirzebruch spectral sequence in oriented bordism starts from the term $E^2_{s,t}=H_s(X,\Omega_t^{SO})$ and converges to the oriented bordism group~$MSO_*(X)$. (Here $\Omega^{SO}_*$ is the oriented bordism ring of a point.) A homology class $z\in H_n(X,\Z)=E^2_{n,0}$ is realisable if and only if it is a cycle of all differentials of the spectral sequence.
Study of the differentials of the Atiyah--Hirzebruch spectral sequence allowed Buchstaber to obtain the following result.

\begin{theorem}[Buchstaber~\cite{Buc69}, 1969]\label{theorem_Buch}
1) Consider $z\in H_n(X,\Z)$. Suppose the group $H_{n-2q(p-1)-1}(X,\Z)$ has no $p$-torsion for every odd prime~$p$ and every~$q\ge 1$ such that $2q(p-1)\ge m$. Then the homology class $z$ is realisable with multiplicity 
$$
\lambda(m)=\prod_{p\ \text{odd prime}}p^{\left[
\frac{m-1}{2(p-1)}
\right]}
$$

2) This estimate is sharp in the sense that for every $m$, there exist a topological space~$X$ and a homology class~$z$ which satisfy the above conditions, but the class~$kz$ is non-realisable for every $k<\lambda(m)$.

3) Every homology class $z\in H_n(X,\Z)$ is realisable with multiplicity~$\lambda(n-1)$.\end{theorem}

The third claim of this theorem gives the best known estimate for the coefficient~$k(n)$ in Theorem~\ref{theorem_Thom}. 

Sullivan~\cite{Sul71} found a geometric interpretation for the differentials of the Atiyah--Hirzebruch spectral sequence. 
Recall that an $n$-dimensional \textit{pseudo-manifold} is a simplicial complex~$Z^n$ such that each simplex of~$Z^n$ is contained in an $n$-simplex and each
$(n-1)$-simplex of~$Z^n$ is contained in exactly two $n$-simplices. All pseudo-manifolds under consideration will be compact. The following proposition easily follows from the definition of singular homology.  

\begin{proposition}\label{propos_pseudo}
For any homology class $z\in H_n(X,\Z)$, there exist an oriented $n$-di\-mensional pseudo-manifold~$Z^n$ and a mapping $f\colon Z^n\to X$ such that $f_*[Z^n]=z$.
\end{proposition}

The approach of Sullivan is to study obstructions to resolving singularities of the pseudo-manifold~$Z^n$.
Suppose that $Z^n$ is a ``manifold up to codimension $s-1$''. This means that the link of every simplex in~$Z^n$ of codimension less or equal to~$s-1$ is PL homeomorphic to the standard sphere of appropriate dimension. Then Sullivan defined a simplicial cycle
$$
\mathfrak{o}_s=\sum_{\codim\sigma=s}[\link\sigma]\sigma\in C_{n-s}(Z^n,\Omega_{s-1}^{SPL})
$$
and proved that its homology class is a complete obstruction to resolving singularities of~$Z^n$, that is, to finding a degree~$1$ mapping $Z_1^n\to Z^n$ such that $Z^n_1$ is a ``manifold up to codimension~$s$''. (Here $\Omega_*^{SPL}$ is a PL oriented bordism ring of a point and square brackets denote the bordism class.) Sullivan showed that in the  situation described, the fundamental homology class~$[Z^n]$ is a cycle of the differentials $d^2,\ldots,d^{s-1}$ of the Atiyah--Hirzebruch spectral sequence, and $d^s[Z^n]$ is equal to the homology class of the cycle~$\mathfrak{o}_s$. He also announced a smooth version of this result. However, an accurate exposition has never appeared.
Notice that Sullivan's result does not allow to prove Theorem~\ref{theorem_Thom} combinatorially without usage of algebraic topology. Actually, it just shifts all difficulties to the computation of the bordism ring of a point. 

\section{Combinatorial realisation of cycles}

In~\cite{Gai07} the author found a purely combinatorial proof of the fact that every integral homology class becomes realisable after multiplication by certain positive integer. This proof provided a direct combinatorial construction that, given a singular simplicial cycle representing a homology class, yields a smooth manifold and a mapping of this manifold that realises the given homology class with certain multiplicity. This approach was further developed in the author's papers~\cite{Gai08a}--\cite{Gai12}. 
Unfortunately, this combinatorial approach does not allow us to obtain a reasonable estimate for the multiplicity with which a homology class is realised.

By Proposition~\ref{propos_pseudo}, every integral homology class $z\in H_n(X,\Z)$ is an image of the fundamental class of an oriented pseudo-manifold~$Z^n$. Moreover, suppose $X$ is arcwise connected, then the pseudo-manifold~$Z^n$ can be chosen \textit{strongly connected\/}, which means that any two $n$-simplices of~$Z^n$ can be connected by a sequence of $n$-simplices so that any two consecutive simplices in this sequence possess a common facet. Our goal is to construct a smooth oriented manifold~$M^n$ and a continuous mapping $f\colon M^n\to Z^n$ such that $f_*[M^n]=q[Z^n]$ for certain $q>0$.

\begin{figure}
\begin{center}
\unitlength=.73mm
\begin{picture}(164,36)

\put(0,0){%
\begin{picture}(0,0)

\put(0,18){\line(1,0){48}}
\put(12,0){\line(1,0){24}}
\put(12,36){\line(1,0){24}}

\put(0,18){\line(2,3){12}}
\put(12,0){\line(2,3){24}}
\put(36,0){\line(2,3){12}}

\put(0,18){\line(2,-3){12}}
\put(12,36){\line(2,-3){24}}
\put(36,36){\line(2,-3){12}}

\end{picture}%
}

\put(64,0){%
\begin{picture}(0,0)

\put(8,18){\line(1,0){8}}
\put(32,18){\line(1,0){8}}
\put(20,0){\line(1,0){8}}
\put(20,36){\line(1,0){8}}

\put(4,24){\line(2,3){4}}
\put(16,6){\line(2,3){4}}
\put(28,24){\line(2,3){4}}
\put(40,6){\line(2,3){4}}

\put(4,12){\line(2,-3){4}}
\put(16,30){\line(2,-3){4}}
\put(28,12){\line(2,-3){4}}
\put(40,30){\line(2,-3){4}}

\thicklines

\put(8,6){\line(1,0){8}}
\put(32,6){\line(1,0){8}}
\put(8,30){\line(1,0){8}}
\put(32,30){\line(1,0){8}}
\put(20,12){\line(1,0){8}}
\put(20,24){\line(1,0){8}}

\put(4,12){\line(2,3){4}}
\put(16,30){\line(2,3){4}}
\put(28,12){\line(2,3){4}}
\put(16,18){\line(2,3){4}}
\put(40,18){\line(2,3){4}}
\put(28,0){\line(2,3){4}}

\put(4,24){\line(2,-3){4}}
\put(16,6){\line(2,-3){4}}
\put(28,24){\line(2,-3){4}}
\put(16,18){\line(2,-3){4}}
\put(40,18){\line(2,-3){4}}
\put(28,36){\line(2,-3){4}}

\end{picture}%
}

\put(120,0){%
\begin{picture}(0,0)

\thicklines

\put(7.67,18.5){\line(1,0){8.66}}
\put(31.67,18.5){\line(1,0){8.66}}
\put(7.67,17.5){\line(1,0){8.66}}
\put(31.67,17.5){\line(1,0){8.66}}
\put(19.67,0.5){\line(1,0){8.66}}
\put(19.67,35.5){\line(1,0){8.66}}

\put(4.33,23.5){\line(2,3){4.33}}
\put(16.33,5.5){\line(2,3){4.33}}
\put(28.33,23.5){\line(2,3){4.33}}
\put(39.33,6){\line(2,3){4.33}}
\put(15.33,6){\line(2,3){4.33}}
\put(27.33,24){\line(2,3){4.33}}

\put(4.33,12.5){\line(2,-3){4.33}}
\put(16.33,30.5){\line(2,-3){4.33}}
\put(28.33,12.5){\line(2,-3){4.33}}
\put(39.33,30){\line(2,-3){4.33}}
\put(15.33,30){\line(2,-3){4.33}}
\put(27.33,12){\line(2,-3){4.33}}
\thicklines

\put(8.66,6){\line(1,0){6.67}}
\put(32.66,6){\line(1,0){6.67}}
\put(8.66,30){\line(1,0){6.67}}
\put(32.66,30){\line(1,0){6.67}}
\put(20.66,12){\line(1,0){6.67}}
\put(20.66,24){\line(1,0){6.67}}

\put(4.33,12.5){\line(2,3){3.33}}
\put(16.33,30.5){\line(2,3){3.33}}
\put(28.33,12.5){\line(2,3){3.33}}
\put(16.33,18.5){\line(2,3){3.33}}
\put(40.33,18.5){\line(2,3){3.33}}
\put(28.33,0.5){\line(2,3){3.33}}

\put(4.33,23.5){\line(2,-3){3.33}}
\put(16.33,5.5){\line(2,-3){3.33}}
\put(28.33,23.5){\line(2,-3){3.33}}
\put(16.33,17.5){\line(2,-3){3.33}}
\put(40.33,17.5){\line(2,-3){3.33}}
\put(28.33,35.5){\line(2,-3){3.33}}

\end{picture}%
}
\put(54,18){\vector(1,0){10}}
\put(112,18){\vector(1,0){10}}

\end{picture}
\end{center}
\caption{Truncating simplices of~$Z^n$ ($n=2$)}\label{fig_trunc_Z}
\end{figure}
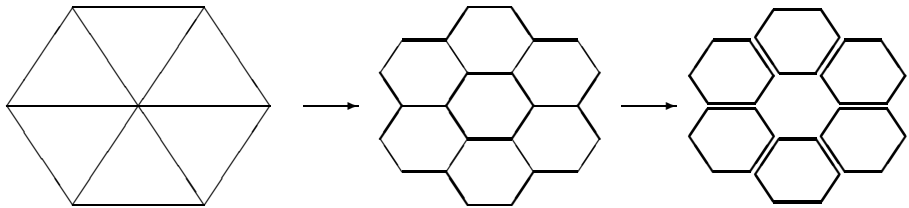

In this section we describe informally main ideas of our construction. The first step of the construction was partially inspired by Sullivan's approach based on resolving singularities of~$Z^n$. Indeed, $Z^n$ can be regarded as  a manifold with singularities in codimension~$2$ skeleton. The worst singularities occur in the vertices of~$Z^n$. Truncate the vertices of~$Z^n$.
This means that we delete from~$Z^n$ a small neibourhood of each vertex such that in each simplex of~$Z^n$ this neighbourhood is bounded by a hyperplane. However, we still have singularities in edges of~$Z^n$. Then we truncate all edges of~$Z^n$. Further, we truncate consecutively all simplices of~$Z^n$ of dimensions~$2,3,\ldots,n-1$, see Fig.~\ref{fig_trunc_Z}. (Though we have no singularities in the interiors of $(n-1)$-simplices of~$Z^n$, it is convenient to truncate them too.) 

This procedure yields a disjoint union of convex polytopes~$\Pi^n_{\sigma}$ corresponding to $n$-simplices~$\sigma$ of~$Z^n$. Each facet of~$\Pi^n_{\sigma}$ corresponds to those subsimplex~$\tau\subset\sigma$ after truncation of which this facet appears. The facet of~$\Pi^n_{\sigma}$ corresponding to~$\tau$ will be denoted by~$F_{\sigma,\tau}$.
Next, for each simplex~$\tau$, we want to pair off and to glue to each other the facets~$F_{\sigma,\tau}$ of all polytopes~$\Pi^n_{\sigma}$, $\sigma\supset\tau$, so that to obtain a smooth manifold without boundary. 
Indeed, we shall replace each~$\Pi^n_{\sigma}$ with several copies of it before such pairings off and gluings, because we are interested in the realisation of a \textit{multiple\/} of the given homology class.

To understand the structure of each  polytope~$\Pi^n_{\sigma}$ we study what happens to an $n$-simplex of~$Z^n$ under the described truncations. It is convenient to identify this simplex with the standard simplex $\Delta^n\subset\R^{n+1}$ given by 
$$
x_1+\cdots+x_{n+1}=1,\quad x_i\ge 0,\ i=1,\ldots,n+1.
$$ 
First, truncate the vertices of~$\Delta^n$ by the hyperplanes $x_i=1-\varepsilon_1$. Second, truncate the edges of~$\Delta^n$ by the hyperplanes $x_{i_1}+x_{i_2}=1-\varepsilon_2$. On the $k$th step, truncate all $(k-1)$-faces of~$\Delta^n$ by the hyperplanes
\begin{equation}\label{eq_Fomega}
x_{i_1}+\cdots+x_{i_k}=1-\varepsilon_k
\end{equation}
The numbers~$\varepsilon_k$ are chosen so that $\varepsilon_1<\frac12$ and
$2\varepsilon_k<\varepsilon_{k-1}$.
It is easy to see that after all such truncations we shall obtain the convex polytope with $(n+1)!$ vertices obtained by permutations of coordinates of the point $$(1-\varepsilon_1,\varepsilon_1-\varepsilon_2,\ldots,\varepsilon_{n-1}-\varepsilon_n,\varepsilon_n).$$
This polytope is called a \textit{permutahedron\/} and is denoted by~$\Pi^n$. Usually, under a permutahedron they mean the convex hull of the $(n+1)!$ points obtained by permutations of coordinates of the point
$(1,2,\ldots,n+1)$. Our permutahedron~$\Pi^n$ is combinatorially equivalent to this standard permutahedron. 
 The polytope~$\Pi^n$ is \textit{simple\/}, i.\,e., every vertex of~$\Pi^n$ is contained in exactly $n$ facets. 
 
 We denote the set $\{1,2,\ldots,n+1\}$ by~$[n+1]$.
 The facets of~$\Pi^n$ are in one-to-one correspondence with non-empty proper subsets $\omega\subset[n+1]$. The facet~$F_{\omega}$ corresponding to a subset $\omega=\{i_1,\ldots,i_k\}$ is given by equation~\eqref{eq_Fomega}. It is convenient to denote by~$\Delta_{\omega}$ the face of~$\Delta^n$ given by the equations $x_i=0$, $i\notin\omega$. Then $\dim\Delta_{\omega}=|\omega|-1$, and the facet~$F_{\omega}$ appears after the truncation of the face~$\Delta_{\omega}$. 

\begin{proposition}\label{propos_mapPiDelta}
There is a piecewise linear mapping $\pi\colon\Pi^n\to\Delta^n$ such that $\pi(F_{\omega})=\Delta_{\omega}$ for every~$\omega$.
\end{proposition}

Thus we have come to the main idea of the construction. For certain positive integer~$q$, we want to replace each $n$-simplex~$\sigma$ of~$Z^n$ with the $q$ copies of the permutahedron obtained from~$\sigma$ by truncating its faces, and then to pair off and to glue to each other the facets of all these permutahedra in a proper way.

\section{Construction}
\label{section_constr}

In this section we shall describe the construction of combinatorial realisation of cycles in details. We shall follow the group-theoretic point of view suggested by the author in~\cite{Gai12}. 
For a given strongly connected oriented pseudo-manifold~$Z^n$, our goal is to construct a smooth manifold~$M^n$ glued from permutahedra and a mapping $f\colon M^n\to Z^n$ such that $f_*[M^n]=q[Z^n]$ for certain positive integer~$q$.

A colouring of vertices of~$Z^n$ in $n+1$ colours  will be called \textit{regular\/} if the vertices of every $n$-simplex are coloured in $n+1$ pairwise distinct colours. We shall conveniently denote the colours by $1,2,\ldots,n+1$. For a subset $\omega\subset[n+1]$, we shall say that a simplex $\tau$ has \textit{type\/}~$\omega$ if the set of colours of vertices of~$\tau$ coincides with~$\omega$.  We shall easily achieve that $Z^n$ has a regular colouring of vertices by replacing~$Z^n$ with its first barycentric subdivision. 
Such colouring of vertices will allow us to identify canonically every simplex~$\sigma$ of~$Z^n$ with the simplex~$\Delta^n$, and every polytope~$\Pi^n_{\sigma}$ with the standard permutahedron~$\Pi^n$. In the sequel we shall fix a regular colouring of verices of~$Z^n$.

We denote by $A$ the set of all $n$-simplices of~$Z^n$.
For every simplex $\sigma\in A$, the colouring of its vertices in colours $1,\ldots,n+1$ provides the orientation of~$\sigma$. On the other hand, the global orientation of~$Z^n$ also provides the orientation of~$\sigma$. We denote by $A_+$ the set of all simplices~$\sigma\in A$ for which the global orientation coincides with the orientation induced by the colouring and denote by~$A_-$ the set of all simplices~$\sigma\in A$ for which these two orientations are opposite to each other.

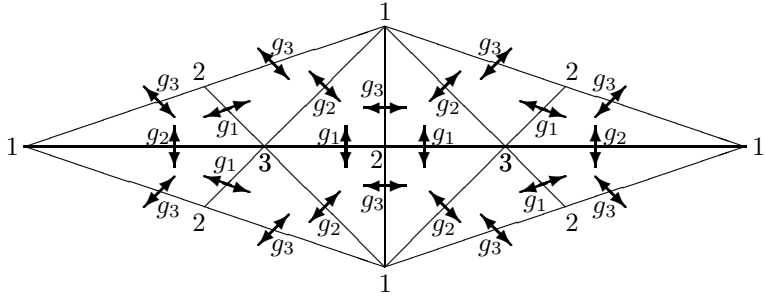
\begin{figure}
\begin{center}
\unitlength=.8mm
\begin{picture}(120,44)

\put(0,20){\line(1,0){120}}
\put(0,20){\line(3,1){60}}
\put(0,20){\line(3,-1){60}}
\put(60,0){\line(3,1){60}}
\put(60,40){\line(3,-1){60}}
\put(60,0){\line(0,1){40}}
\put(60,0){\line(1,1){30}}
\put(60,0){\line(-1,1){30}}
\put(60,40){\line(1,-1){30}}
\put(60,40){\line(-1,-1){30}}

\thicklines

\put(67.5,27.5){\vector(1,1){5}}
\put(72.5,32.5){\vector(-1,-1){5}}

\put(52.5,12.5){\vector(-1,-1){5}}
\put(47.5,7.5){\vector(1,1){5}}

\put(67.5,12.5){\vector(1,-1){5}}
\put(72.5,7.5){\vector(-1,1){5}}

\put(52.5,27.5){\vector(-1,1){5}}
\put(47.5,32.5){\vector(1,-1){5}}

\put(63.5,26.5){\vector(-1,0){7}}
\put(56.5,26.5){\vector(1,0){7}}

\put(63.5,13.5){\vector(-1,0){7}}
\put(56.5,13.5){\vector(1,0){7}}

\put(66.5,16.5){\vector(0,1){7}}
\put(66.5,23.5){\vector(0,-1){7}}

\put(53.5,16.5){\vector(0,1){7}}
\put(53.5,23.5){\vector(0,-1){7}}

\put(95,16.5){\vector(0,1){7}}
\put(95,23.5){\vector(0,-1){7}}

\put(25,16.5){\vector(0,1){7}}
\put(25,23.5){\vector(0,-1){7}}

\put(82.5,12.5){\vector(3,1){7.5}}
\put(90,15){\vector(-3,-1){7.5}}

\put(82.5,27.5){\vector(3,-1){7.5}}
\put(90,25){\vector(-3,1){7.5}}

\put(37.5,12.5){\vector(-3,1){7.5}}
\put(30,15){\vector(3,-1){7.5}}

\put(37.5,27.5){\vector(-3,-1){7.5}}
\put(30,25){\vector(3,1){7.5}}

\put(20,10){\vector(1,1){5}}
\put(25,15){\vector(-1,-1){5}}
\put(39,3.7){\vector(1,1){5}}
\put(44,8.7){\vector(-1,-1){5}}

\put(100,10){\vector(-1,1){5}}
\put(95,15){\vector(1,-1){5}}
\put(81,3.7){\vector(-1,1){5}}
\put(76,8.7){\vector(1,-1){5}}

\put(20,30){\vector(1,-1){5}}
\put(25,25){\vector(-1,1){5}}
\put(39,36.3){\vector(1,-1){5}}
\put(44,31.3){\vector(-1,1){5}}

\put(100,30){\vector(-1,-1){5}}
\put(95,25){\vector(1,1){5}}
\put(81,36.3){\vector(-1,-1){5}}
\put(76,31.3){\vector(1,1){5}}

\put(22,9){$g_3$}
\put(41,2.7){$g_3$}
\put(94.5,9){$g_3$}
\put(75.5,2.7){$g_3$}

\put(22,30){$g_3$}
\put(41,36){$g_3$}
\put(94.5,30){$g_3$}
\put(75.5,36){$g_3$}

\put(-3,18.5){$1$}
\put(121,18.5){$1$}
\put(59,-4.2){$1$}
\put(59,41.1){$1$}

\put(28,6){$2$}
\put(90,6){$2$}
\put(28,31){$2$}
\put(90,31){$2$}
\put(57.5,16.5){$2$}

\put(39,16){$3$}
\put(39,16){$3$}

\put(79,16){$3$}
\put(79,16){$3$}

\put(48.8,21){$g_1$}
\put(67.8,21){$g_1$}

\put(20.3,21){$g_2$}
\put(96.3,21){$g_2$}

\put(56,29){$g_3$}
\put(56,10){$g_3$}

\put(31.5,16.5){$g_1$}
\put(83,9.9){$g_1$}

\put(32,22.5){$g_1$}
\put(85,22.5){$g_1$}

\put(48.8,6){$g_2$}
\put(67.8,6){$g_2$}

\put(48,26){$g_2$}
\put(68.5,26){$g_2$}

\end{picture}
\end{center}

\caption{Action of~$G$ on~$A$}\label{fig_lambda}
\end{figure}

 Let $G$ be the free product of $n+1$ copies of the group~$\Z_2$, that is,
$$
G=\langle
g_1,\ldots,g_{n+1}\mid g_1^2=\cdots=g_{n+1}^2=1
\rangle.
$$ 
Consider the action of~$G$ on $A$ such that $g_i\sigma$ is a unique simplex in~$A$ that is distinct from~$\sigma$ and has a common facet of type~$[n+1]\setminus\{i\}$ with $\sigma$, see Fig.~\ref{fig_lambda}. Obviously, $g_i\sigma\in A_-$ whenever $\sigma\in A_+$, and vice versa. 
Choose an arbitrary simplex $\sigma_0\in A_+$ and let $H\subset G$ be the stabilizer of~$\sigma_0$. We have the natural identification $A=G/H$. Moreover,  the pseudo-manifold~$Z^n$ is given by
$$
Z^n=(\Delta^n\times(G/H))/\sim,
$$
where $\sim$ is the equivalence relation such that $(p,gH)\sim (p',g'H)$ if and only if $p=p'$ and $g'H=ygH$ for  an element~$y$ of the subgroup of~$G$ generated by all~$g_i$ such that  $x\in\Delta_{[n+1]\setminus\{i\}}$. The point of~$Z^n$ corresponding to the equivalence class of a pair $(p,gH)$ will be denoted by~$[p,gH]$.

Let $\Omega$ be  the set of all non-empty proper subsets of~$[n+1]$ partially ordered by inclusion. The cardinality of~$\Omega$ is equal to $m=2^{n+1}-2$.  Let $\F$ be the free product of $m$ copies of the group~$\Z_2$
labeled by the elements of the set~$\Omega$. The generator of the factor~$\Z_2$ corresponding to an element~$\omega\in\Omega$ will be denoted by~$x_{\omega}$. Thus,
$$
\F=\langle x_{\omega},\omega\in\Omega\mid x_{\omega}^2=1\rangle.
$$
For each $\omega\in\Omega$, we consider the automorphism $\psi_{\omega}\in\Aut(\F)$ such that
$$
\psi_{\omega}(x_{\gamma})=\left\{
\begin{aligned}
x_{\omega}&x_{\gamma}x_{\omega}&&\text{if}\ \gamma\subset\omega,\\
&x_{\gamma}&&\text{if}\ \gamma\not\subset\omega
\end{aligned}
\right.
$$ 
for every $\gamma\in\Omega$. We denote by $\Psi$ the subgroup of~$\Aut(\F)$ generated by the elements~$\psi_{\omega}$, $\omega\in\Omega$.
Now, we consider the semi-direct product $\F\rtimes\Psi$ corresponding to the tautological action of~$\Psi\subset\Aut(\F)$ on~$\F$ and we consider the elements $s_{\omega}=x_{\omega}\psi_{\omega}\in \F\rtimes\Psi$. Denote by $S$ the set consisting of $m$ elements~$s_{\omega}$ and denote by~$W$ the subgroup of~$\F\rtimes\Psi$ generated by the elements~$s_{\omega}$. 
\begin{proposition}[\cite{Gai12}]\label{propos_Coxeter}
The elements $s_{\omega}$ satisfy the relations $s_{\omega}^2=1$ for every $\omega\in\Omega$ and $s_{\omega_1}s_{\omega_2}=s_{\omega_2}s_{\omega_1}$ whenever $\omega_1\subset\omega_2$. These relations imply all relations among the elements $s_{\omega}$. Thus $(W,S)$ is a right-angular Coxeter system. 
\end{proposition}

Davis~\cite{Dav83} assigned a cell complex to every Coxeter system. For the Coxeter system~$(W,S)$ this construction yields the manifold
$$
\widetilde{M}^n=(\Pi^n\times W)/\sim,
$$ 
where $\sim$ is the equivalence relation such that $(p,w)\sim(p',w')$ if and only if $p=p'$ and  $w'w^{-1}$ belongs to the subgroup of~$W$ generated by all $s_{\omega}$ such that $p\in F_{\omega}$. 
The point of~$\widetilde{M}^n$ corresponding to the equivalence class of a pair $(p,w)$ will be denoted by~$[p,w]$. 
The right action of~$W$  on~$\widetilde{M}^n$ is given by~$[p,w]\cdot w'=[p,ww']$.
Davis proved that the manifold~$\widetilde{M}^n$ admits a canonical smoothing and is contractible.

For each non-empty proper subset $\omega\subset[n+1]$, we denote by $\F_{\omega}$ the subgroup of~$\F$ generated by all elements $x_{\omega'}$ such that $\omega'\supset\omega$ and we denote by $G_{\omega}$ the subgroup of~$G$ generated by all elements~$g_i$ such that $i\in\omega$.

Now, we consider the mapping $\F\rtimes\Psi\to \F$ such that $\psi x\mapsto x$ for every~$x\in \F$ and~$\psi\in\Psi$ and we denote by~$\theta$ the restriction of this mapping to~$W$. The mapping~$\theta$ is not a homomorphism. However, it satisfies the following property. 

\begin{proposition}\label{propos_theta}
Suppose $w\in W$ and $\omega\in\Omega$; then there is an element $y\in \F_{\omega}$ such that
$
\theta(s_{\omega}w)=yx_{\omega}y^{-1}\theta(w).
$
\end{proposition} 

Consider the homomorphism $\lambda\colon \F\to G$ given by $\lambda(x_{\omega})=g_{\min([n+1]\setminus\omega)}$, where $\min(\gamma)$ is the minimal element in~$\gamma$. (The minimum is not important here. Indeed, we could take for~$\lambda(x_{\omega})$ an arbitrary element~$g_i$ such that $i\notin\omega$.) Let $\varphi=\lambda\theta$. Proposition~\ref{propos_theta} immediately implies

\begin{corollary}\label{cor_phi}
Suppose $w\in W$ and $\omega\in\Omega$; then there is an element $g\in G_{[n+1]\setminus\omega}$ such that
$
\varphi(s_{\omega}w)=g\varphi(w).
$
Besides, $g$ is a product of odd number of generators~$g_i$.
\end{corollary}

\begin{corollary}\label{cor_phi_action}
Suppose $w\in W$ and $\omega\in\Omega$; then the simplices $\sigma=\varphi(w)\sigma_0$ and $\sigma'=\varphi(s_{\omega}w)\sigma_0$ have a common face of type~$\omega$. Besides, $\sigma'\in A_-$ if $\sigma\in A_+$ and $\sigma'\in A_+$ if $\sigma\in A_-$.
\end{corollary}

Let $\pi\colon\Pi^n\to \Delta^n$ be the mapping in Proposition~\ref{propos_mapPiDelta}. Define the mapping $\Phi\colon\widetilde{M}^n\to Z^n$ by
$$
\Phi([p,w])=[\pi(p),\varphi(w)].
$$
Corollary~\ref{cor_phi_action} implies that $\Phi$ is well defined and maps every cell of~$\widetilde{M}^n$ onto a simplex of~$Z^n$ preserving the orientation. However, the manifold~$\widetilde{M}^n$ is not compact. Now we need to find a subgroup~$\Gamma\subset W$ of finite index such that $\Gamma$ acts freely on~$\widetilde{M}^n$ and the mapping~$\Phi$ is $\Gamma$-invariant, i.\,e., passes through $M^n=\widetilde{M}^n/\Gamma$:
$$
\Phi\colon\widetilde{M}^n\to\widetilde{M}^n/\Gamma\xrightarrow{f}Z^n.
$$
Then we shall see that $f_*[M^n]=q[Z^n]$ for a positive integer~$q$.

The subgroup~$\Gamma$ is constructed as follows.
First, consider the group $\overline{H}=\lambda^{-1}(H)$. Obviously, $\overline{H}$ has finite index in~$\F$. Second, consider the subgroup $\Psi_H\subset\Psi$ consisting of all automorphisms~$\psi$ such that $\psi(x)\overline{H}=x\overline{H}$ for every~$x\in\F$. 
It is easy to show that the subgroup~$\Psi_H$ has finite index in~$\Psi$. Besides, we see that $\psi(\overline{H})=\overline{H}$ for every $\psi\in\Psi_H$. Hence, $\overline{H}\rtimes\Psi_H$ is a well-defined subgroup of finite index in~$\F\rtimes\Psi$. We put $W_H=W\cap(\overline{H}\rtimes\Psi_H)$. Then $W_H$ is a finite index subgroup of~$W$.

\begin{proposition}
The mapping $\Phi$ is $W_H$-invariant. 
\end{proposition}

Now consider the group $\Z_2^n=\Z_2\times\cdots\times\Z_2$ with generators $a_1,\ldots,a_n$ and define the homomorphism $\eta\colon W\to\Z_2^n$ by $\eta(s_{\omega})=a_{|\omega|}$. The kernel $K=\ker\eta$ acts freely on~$\widetilde{M}^n$. Therefore, $\Gamma=W_H\cap K$ is the required subgroup. Put $M^n=\widetilde{M}^n/\Gamma$ and let $f\colon M^n\to Z^n$ be the quotient of~$\Phi$. Then $f_*[M^n]=q[Z^n]$, where $q=\frac{|W:\Gamma|}{|G:H|}$. 

\section{Tomei manifold}

Let us make the following important observation. Consider the manifold $M^n_0=\widetilde{M}^n/K$. Since the action of~$K$ on~$\widetilde{M}^n$ is free and preserves the orientation, we see that $M^n_0$ is an oriented smooth closed manifold glued out of $2^n$ permutahedra. Notice that for each oriented pseudo-manifold~$Z^n$, the manifold $M^n=\widetilde{M}^n/\Gamma$ obtained by the construction described in the previous section is a finite-sheeted covering of~$M^n_0$. Thus we obtain the following 

\begin{theorem}[\cite{Gai08a}, \cite{Gai08b}]\label{theorem_M0}
For each topological space~$X$ and each homology class $z\in H_n(X,\Z)$, there are a finite-sheeted covering~$M^n$ of~$M^n_0$ and a continuous mapping $f\colon M^n\to X$ such that $f_*[M^n]=qz$ for a positive integer~$q$. 
\end{theorem} 

The manifold~$M_0^n$ has many interesting properties. Tomei~\cite{Tom84} proved that $M_0^n$ is homeomorphic to the isospectral manifold of symmetric tridiagonal real $(n+1)\times(n+1)$-matrices, that is, to the manifold consisting of all symmetric tridiagonal $(n+1)\times(n+1)$-matrices
$$
\left(
\begin{array}{ccccc}
a_1&b_1&0&\cdots&0\\
b_1&a_2&b_2&\cdots&0\\
0&b_2&a_3&\cdots&0\\
\vdots&\vdots&\vdots&\ddots&\vdots\\
0&0&0&\cdots&a_{n+1}
\end{array}
\right),\qquad
a_i,b_i\in\R,
$$ 
with the fixed simple spectrum $\lambda_1<\lambda_2<\cdots<\lambda_{n+1}$. Hence $M^n_0$ will be called the \textit{Tomei manifold\/}.
Besides, $M^n_0$ is the so-called \textit{small cover induced from the linear model\/} of the permutahedron~$\Pi^n$ (see~\cite{DaJa91a}).
The result of Davis~\cite{Dav83} claims that the manifold~$\widetilde{M}^n$ is contractible. Hence, the manifold~$M^n_0$ is \textit{aspherical\/}, i.e., $\pi_i(M^n_0)=0$ for $i>1$. 

\section{Domination relation and URC-manifolds}

Let $M^n$ and $N^n$ be oriented closed manifolds of the same dimension~$n$. We say that $M^n$ \textit{dominates\/}~$N^n$ and write $M^n\geqslant N^n$ if there exists a non-zero degree mapping $M^n\to N^n$. Domination is a transitive relation on the set of homotopy types of oriented connected closed manifolds. This relation goes back to  the works of Milnor and Thurston~\cite{MiTh77} and Gromov~\cite{Gro82} and was first explicitly defined in the paper of Carlson and Toledo~\cite{CaTo89} with a reference to a lecture of Gromov.  Obviously, $M^n\geqslant S^n$ for every $n$-dimensional manifold~$M^n$. Hence the sphere~$S^n$ is the minimal element with respect to the domination relation.   
We shall focus on the following question.

\begin{question}[Carlson--Toledo~\cite{CaTo89}, 1989]
\label{quest_main} 
Is there an easily describable maximal class of homotopy types with respect to the domination relation, i.\,e., a collection~$\CC_n$ of $n$-dimensional manifolds such that given any manifold~$N^n$, there is an $M^n\in\CC_n$ satisfying $M^n\geqslant N^n$?
\end{question}

This question is closely related to the problem of finding a collection of manifolds sufficient for  realisation with multiplicities of all integral homology classes. Suppose we have a class~$\CC_n$ satisfying the requirements of Question~\ref{quest_main}. Thom's theorem~\cite{Tho54} (see Theorem~\ref{theorem_Thom}) claims that  for each integral homology class~$z$, a multiple~$kz$ can be realised by an image of the fundamental class of an oriented manifold~$N^n$. Then there exists a manifold $M^n\in \CC_n$ such that $M^n\geqslant N^n$. Therefore, a multiple of~$z$ is realisable by an image of the fundamental class of~$M^n$.

Let us mention some known results about the domination relation. In dimension~$2$ the structure of the domination relation is clear. In dimension~$3$, there are a lot of results.   
We do not want to describe them here. A good survey can be found in~\cite{Wan02}.
However, not so much is known in dimensions $n\ge 4$. 

First, let us mention some negative results. It is known that there exist manifolds that cannot be dominated by a K\"ahler manifold~\cite{CaTo89} and there exist manifolds that cannot be dominated by a direct product $M_1\times M_2$ of two manifolds of positive dimensions~\cite{KoLo08}. A survey of other results in this direction can be found in~\cite{KoLo08}.

For a long time, positive results  were related to the so-called \textit{hyperbolisation procedure\/} introduced by 
Gromov~\cite{Gro87}. For each polyhedron~$P$, this procedure yields a polyhedron~$P_h$ and a mapping~$f\colon P_h\to P$ such that $P_h$ has a metric of non-positive curvature in sense of Alexandrov, $f$ induces an injection in integral cohomology, and singularities of~$P_h$ are ``not worse'' than  singularities of~$P$. In particular, for an oriented PL manifold~$N^n$, this construction yields a PL manifold~$N^n_h$ and a degree~$1$ mapping $f\colon N^n_h\to N^n$. Besides,  $f^*$ takes the rational Pontryagin classes of~$N^n$ to the rational Pontryagin classes of~$N^n_h$. This result implies that every manifold can be dominated with degree~$1$ by a manifold of non-positive curvature in sense of Alexandrov. Gromov's hyperbolisation procedure has been improved by Davis--Januszkiewicz~\cite{DaJa91b}, Charney--Davis~\cite{ChDa95}, and other authors. Recently Ontaneda~\cite{Ont11} has proved that, for a smooth manifold~$N^n$, its hyperbolisation~$N^n_h$ can be constructed to be a Riemannian manifold of negative sectional curvature in an arbitrarily small interval~$[-1-\varepsilon,-1]$.

Thus the class of all $n$-dimensional manifolds of sectional curvature in the interval $[-1-\varepsilon,-1]$ can be taken for~$\CC_n$ in the Carlson--Toledo question. However, this class still seems to be too wide for Question~\ref{quest_main}. Actually, there is an essential difference between hyperbolisation and domination. The condition that 
$f\colon M^n\to N^n$ is a domination is certainly much weaker than the condition that $f$ is a degree~$1$ map that preserves the rational Pontryagin classes. Hence it is natural to hope that we would be able to find a smaller class~$\CC_n$. However, such a class had not been known until the author~\cite{Gai08a} proved Theorem~\ref{theorem_M0}. This theorem easily implies that for~$\CC_n$ one can take the class of all finite-sheeted coverings of the Tomei manifold~$M_0^n$. It is convenient to formulate this result in terms of a so-called \textit{virtual domination relation\/}. Let $M^n$ and $N^n$ be oriented closed manifolds of the same dimension. We shall say that $M^n$ \textit{virtually dominates\/}~$N^n$ if a finite-sheeted covering of~$M^n$ dominates~$N^n$. The virtual domination relation is a transitive relation on the set of homotopy types of oriented closed $n$-dimensional manifolds. 

\begin{corollary}
The Tomei manifold~$M_0^n$ is maximal with respect to the virtual domination relation, i.\,e., $M_0^n$ virtually dominates every oriented closed $n$-dimensi\-onal manifold.
\end{corollary}

In~\cite{Gai12} the author suggested to study the class of manifolds that are maximal with respect to the virtually domination relation. Such manfolds~$M^n$ will be called \textit{URC-manifolds\/} because they satisfy the following

\begin{URC} 
For every topological space~$X$ and every homology class $z\in H_n(X,\Z)$, there exist a finite-sheeted covering~$\hM^n$ of~$M^n$ and a mapping $f\colon\hM^n\to X$ such that $f_*\bigl[\hM^n\bigr]=kz$ for a non-zero~$k$.
\end{URC}

\begin{proposition}\label{propos_obvious}
1) Suppose $\hM^n$ is a finite-sheeted covering of~$M^n$; then $\hM^n$ is a URC-manifold if and only if $M^n$ is a URC-manifold.

2) Suppose $M^n$ is a URC-manifold; then the connected sum $M^n\# N^n$ is a URC-manifold for any~$N^n$. 

3) Suppose $N^n$ is a URC-manifold; then $M^n$ is a URC-manifold if and only if $M^n$ virtually dominates~$N^n$. In particular, a manifold~$M^n$ is a URC-manifold if and only if it dominates the Tomei manifold~$M_0^n$.
\end{proposition} 

\begin{corollary}
A class~$\CC_n$ of $n$-dimensional manifolds is maximal with respect to the domination relation, i.\,e., satisfies the requirements of Question~\ref{quest_main} if and only if it contains a URC-manifold. 
\end{corollary}

In~\cite{Gai12} the author found a series of examples of URC-manifolds. These manifolds are obtained by the following general construction due to Davis and Januszki\-ewicz \cite{DaJa91a}. Let $P$ be an $n$-dimensional simple polytope with $m$ facets $F_1,\ldots,F_m$. Consider the group $\Z_2^m=\Z_2\times\cdots\times\Z_2$ with generators $b_1,\ldots,b_m$. 
To the polytope~$P$ corresponds the $n$-dimensional manifold
$$
\RZ_P^n=(P\times\Z_2^m)/\sim,
$$
where $(p,g)\sim (p',g')$ if and only if $p=p'$ and $g'g^{-1}$ belongs to the subgroup of~$\Z_2^m$ generated by all~$b_i$ such that~$x\in F_i$. Since the polytope~$P$ is a smooth manifold with corners, the manifold~$\RZ^n_P$ can be endowed with a smooth structure (see~\cite{Dav83}, \cite{DaJa91a} for details).  

The group~$\Z_2^m$ naturally acts on~$\RZ^n_P$. 
If there exists a subgroup $\Lambda\subset\Z_2^m$ which is isomorphic to~$\Z_2^{m-n}$ and acts freely on~$\RZ^n_P$, then the quotient $\RZ^n_P/\Lambda$ is called a \textit{small cover\/} of~$P$. Small covers of simple polytopes play an important role in \textit{toric topology\/} since thay are real analogues of the so-called \textit{quasitoric manifolds\/}, which are topological analogues of \textit{toric varietes\/} (see~\cite{DaJa91a},~\cite{BuPa02}). Notice that there exist polytopes that admit no small covers.

In~\cite{Gai12} the author has found several sufficient conditions for~$\RZ^n_P$ to be a URC-manifold. To formulate this result we need some definitions and notation. 

By $K_P$ we denote the boundary of the simplicial polytope dual to~$P$. 

A simple polytope~$P$ is called a \textit{flag polytope\/} if $F_{i_1}\cap\ldots\cap F_{i_k}\ne\emptyset$ for any facets $F_{i_1},\ldots,F_{i_k}$ of~$P$ with non-empty pairwise intersections. The dual terminology for the simplicial complex~$K_P$ is that  $K_P$ is a \textit{flag complex\/}, or contains \textit{no empty simplices\/}, or satisfies \textit{no-$\triangle$-condition\/}. We shall say that a flag simple polytope~$P$ contains \textit{no empty $4$-circuit\/} 
if for any its facets $F_{i_1}$, $F_{i_2}$, $F_{i_3}$, and $F_{i_4}$ such that the intersections
$F_{i_1}\cap F_{i_2}$,
$F_{i_2}\cap F_{i_3}$,
$F_{i_3}\cap F_{i_4}$, and
$F_{i_4}\cap F_{i_1}$ are non-empty,
at least one of the two intersections $F_{i_1}\cap F_{i_3}$ and $F_{i_2}\cap F_{i_4}$ is also non-empty.
The dual condition for~$K_P$ is sometimes called \textit{Siebenmann's no-$\square$-condition\/} (see~\cite{Gro87}).

By $\Lob^n$ we denote the $n$-dimensional Lobachevsky space. A convex polytope $P\subset\Lob^n$ is called \textit{right-angular\/} if all its dihedral angles are equal to~$\frac{\pi}{2}$. All polytopes in~$\Lob^n$ under consideration are supposed to be compact. It is well known that a compact right-angular polytope is simple.

\begin{theorem}[\cite{Gai12}]\label{theorem_sc}
Let $P$ be an $n$-dimensional simple polytope  satisfying one of the following conditions:
\begin{enumerate}
\item There exists a simplicial mapping $f\colon K_P\to (\partial \Delta^n)'$ of non-zero degree, where  $(\partial \Delta^n)'$ is the barycentric subdivision of the boundary of the $n$-simplex. (In particular, this holds if~$K_P$ is isomorphic to the barycentric subdivision a simplicial complex.)
\item $P$ is a  simple convex polytope in~$\Lob^n$ such that the interior of~$P$ contains a closed ball of radius 
\begin{equation*}
\rho_n=
\log\left(\sqrt{\frac{n(n+1)(n+2)}{6}}+
\sqrt{\frac{n(n+1)(n+2)}{6}-1}\right).
\end{equation*}
\item $P$ is a compact right-angular convex polytope in~$\Lob^n$.
\item $P$ is a flag simple polytope without an empty $4$-circuit.
\end{enumerate}
Then $\RZ_P^n$ is a URC-manifold. Hence every small cover of~$P$ is a URC-manifold.  
\end{theorem}

This theorem is proved by showing that in each of these cases the manifold~$\RZ^n_P$ virtually dominates the Tomei manifold~$M^n_0$. Notice that in cases~(3) and~(4), the crucial point is that the manifold~$\RZ_P^n$ has a metric of strictly negative curvature in sense of Alexandrov, i.\,e., the universal covering~$\widetilde{\RZ}^n_P$ is a $C\!AT(-1)$-manifold. 
Case~(3) is  related to the following interesting conjecture.

\begin{conjecture}[Kotschick--L\"oh~\cite{KoLo08}, 2008]\label{conj_main}
For $n\ge 2$, every $n$-dimensional oriented closed manifold can be dominated by a hyperbolic manifold, i.\,e., a manifold that admits a Riemannian metric of constant negative curvature. 
\end{conjecture}

This conjecture is trivial for $n=2$ and is known to be true for $n=3$ by a result of Brooks~\cite{Bro85}.
It is easy to show that if~$P\subset \Lob^n$ is a right-angular convex polytope, then $\RZ^n_P$ is a hyperbolic manifold. Actually, in this case $\RZ^n_P$ is isometric to the quotient $\Lob^n/\Gamma$, where $\Gamma$ is a finite index torsion-free subgroup of the uniform trasformation group $W\subset\Isom(\Lob^n)$ generated by reflections in facets of~$P$. 

\begin{corollary}
Suppose $W\subset\Isom(\Lob^n)$ is a uniform right-angular reflection group and $\Gamma$ is a finite index torsion-free subgroup of~$W$. Then the hyperbolic manifold~$\Lob^n/\Gamma$ is a URC-manifold.
\end{corollary}

It is well known that compact right-angular convex polytopes in~$\Lob^n$ exist for~$n=2,3,4$. In particular, in~$\Lob^4$ there exists a regular right-angular $120$-cell. Hence, in dimensions $2$, $3$, and $4$ there are hyperbolic URC-manifolds.

\begin{corollary}
Conjecture~\ref{conj_main} is true in dimension~$4$.
\end{corollary}   

By a result of Vinberg~\cite{Vin84}, for $n\ge 5$, there exist no compact right-angular convex polytopes in~$\Lob^n$, hence, there exist no uniform right-angular reflection subgroups~$W\subset\Isom(\Lob^n)$. So Conjecture~\ref{conj_main} remains open for $n\ge 5$. 

\section{Simplicial volume and URC-manifolds}

Recall that the \textit{simplicial $\ell^1$-semi-norm\/}  in the homology of a topological space~$X$ is defined in the following way. 
Denote by $C_k(X,\R)$ the $k$-dimensional singular simplicial chain group of~$X$ with real coefficients. For each chain~$\zeta=\sum\alpha_i\sigma_i$, where $\alpha_i\in\R$ and $\sigma_i$ are singular simplices, the $\ell^1$-norm of it is equal to
$$
\|\zeta\|_1=\sum|\alpha_i|.
$$ 
The simplicial $\ell^1$-semi-norm of a homology class~$z\in H_k(X,\R)$ is given by
$$
\| z\|_1=\inf\|\zeta\|_1,
$$ 
where the infimum is taken over all singular cycles~$\zeta$ representing the homology class~$z$.
The simplicial $\ell^1$-semi-norm of the fundamental class~$[M^n]$ of an oriented closed manifold~$M^n$ is called the \textit{simplicial volume\/} (or the \textit{Gromov norm\/}) of~$M^n$, and is denoted by~$\|M^n\|$ (see~\cite{Gro82} for details). 

It is well-known that the simplicial $\ell^1$-semi-norm and the simplicial volume have the following properties:

1) If $f\colon X\to Y$ is a continuous mapping and $z\in H_n(X,\R)$, then $\|f_*z\|\le\|z\|$. In particular, if $f\colon M^n\to N^n$ is a mapping of degree~$k$, then $\|M^n\|\ge |k|\|N^n\|$. 

2) If $\widehat{M}^n$ is a $k$-sheeted covering of~$M^n$, then $\|\widehat{M}^n\|=k\|M^n\|$. 

Since in each dimension $n\ge 2$ there exist manifolds of non-zero simplicial volume, these properties imply that $\|M^n\|>0$ for every URC-manifold~$M^n$. The inverse is not true, since, for example, the direct products of surfaces of genera greater or equal to~$2$ have positive simplicial volumes, but are not URC-manifolds.  

Now, for each URC-manifold~$M^n$, we define the corresponding semi-norm \mbox{$\|\cdot\|_{M^n}$} in homology in the following way. For any~$X$ and $z\in H_n(X,\Z)$, we put
$$
\|z\|_{M^n}=\inf\frac{k}{|q|},
$$
where the infimum is taken over all mappings $f\colon\widehat{M}^n\to X$ such that $\hM^n$ is a $k$-sheeted covering of~$M^n$ and $f_*[\hM^n]=qz$, $q\ne 0$. If $X$ is a $3$-manifold and $M^n$ is an oriented surface we (up to a multiplicative constant) obtain the usual Thurston semi-norm.

\begin{theorem}
For each URC-manifold~$M^n$ there exist positive constants $c_1(M^n)$ and $c_2(M^n)$ such that
for every~$X$ and every $z\in H_n(X,\Z)$, we have
\begin{equation}\label{eq_ineq}
c_1(M^n)\|z\|_1\le\|z\|_{M^n}\le c_2(M^n)\|z\|_1.
\end{equation}
\end{theorem}

\noindent{\textit{Proof\/}. Obviously, $\|z\|_1\le \|M^n\|\|z\|_{M^n}$. Hence we can take $c_1(M^n)=\|M^n\|^{-1}$. The second inequality is not so easy. First, let $M^n_0$ be the Tomei manifold. Then $\|z\|_{M^n}\le \|M^n_0\|_{M^n}\|z\|_{M^n_0}$. Hence we suffice to prove the required inequality for the Tomei manifold~$M^n_0$ only.
Second, it is not hard to prove the following}

\begin{proposition}
Let $z\in H_n(X,\Z)$ and let $\varepsilon>0$. Then there exists an oriented $n$-dimensional pseudo-manifold~$Z^n$ with a regular colouring of vertices in $n+1$ colours, and a mapping $g\colon Z^n\to X$ such that $g_*[Z^n]=lz$, $l>0$, and the number of $n$-dimensional simplices of~$Z^n$ is less than $(n+1)!\,l\,(\|z\|_1+\varepsilon)$.
\end{proposition}

The multiple $(n+1)!$ corresponds to the barycentric subdivision which is needed to obtain a regular colouring of vertices.

Apply the construction described in section~\ref{section_constr} to the pseudo-manifold~$Z^n$. We obtain a manifold~$\hM^n_0$ which is a $k$-sheeted covering of~$M_0^n$ and a mapping $f\colon \hM^n_0\to Z^n$ such that $f_*[\hM^n_0]=q[Z^n]$ for a positive~$q$. Suppose, $r$ is the number of $n$-dimensional simplices in~$Z^n$. It immediately follows from the construction that the obtained cell decomposition of~$\hM^n_0$ consists out of~$qr$ permutahedra. This decomposition is the $k$-sheeted covering of the decomposition of~$M_0^n$ into $2^n$ permutahedra. Therefore, $qr=2^nk$. Since $r<(n+1)!\,l\,(\|z\|_1+\varepsilon)$, we see that 
$$
\frac{k}{ql}<\frac{(n+1)!}{2^n}(\|z\|_1+\varepsilon)
$$
But $(gf)_*[\hM^n_0]=qlz$. Hence we obtain that
$
\|z\|_{M^n_0}\le \frac{(n+1)!}{2^n}\|z\|_1
$.
Therefore,
$$
\|z\|_{M^n}\le \frac{(n+1)!}{2^n}\|M^n\|_{M^n_0}\|z\|_1.\eqno{\square}
$$

\section{Open problems}

\begin{question}
Is it possible to modify the construction in section~\ref{section_constr} so that to obtain a construction that realises combinatorially modulo~$2$ homology classes?
\end{question}

A geometric proof of the fact that every modulo~$2$ homology class is realisable was obtained by Buoncristiano and Hacon~\cite{BuHa81}. This proof does not use results of algebraic topology. However, it does not give a combinatorial construction for the manifold realising the homology class. 

Until now our approach to combinatorial realisation of homology classes allows us to say almost nothing about the multiplicities with which the homology classes are realised. So the following question seems to be important.

\begin{question}\label{quest_bound}
Let $M^n$ be a URC-manifold, $n\ge 3$. Does there exist a constant $q(M^n)$ depending only on the manifold~$M^n$ itself such that for every space~$X$ and every homology class~$z\in H_n(X,\Z)$, the class $q(M^n)z$ can be realised by an image of the fundamental class of a finite-sheeted covering of~$M^n$? If ``yes'', does this constant depend on the URC-manifold~$M^n$? What is the minimum of~$q(M^n)$ over all $n$-dimensional URC-manifolds~$M^n$?
\end{question}

Obviously, if such constant~$q(M^n)$ exists for a URC-manifold~$M^n$, then it exists for any URC-manifold of the same dimension.

Since for $n\ge 7$, there are non-realisable homology classes, we certainly cannot achieve $q(M^n)=1$. Nevertheless, it is interesting to understand if it is possible to achieve this replacing the manifold~$M^n$ with a pseudo-manifold.

\begin{question}\label{quest_pseudo}
For $n\ge 3$, does there exist a universal oriented pseudo-manifold~$Z^n$ such that for every~$X$ and every $z\in H_n(X,\Z)$, the class~$z$ can be realised by an image of the fundamental class of a finite-sheeted covering of~$Z^n$?  
\end{question}

A positive answer to Question~\ref{quest_pseudo} will imply a positive answer to Question~\ref{quest_bound}.

\begin{problem}
Characterise simple polytopes~$P$ such that $\RZ^n_P$ is a URC-manifold.
\end{problem}

In~\cite{Gai12} the author has shown that $\RZ^n_P$ is a URC-manifold if it admits a $\Z_2^m$-invariant metric of strictly negative curvature in sense of Alexandrov
(with at least one ``smooth'' point). Here negative curvature cannot be replaced with non-positive curvature. It is well known that the manifold~$\RZ^n_P$ has a $\Z_2^m$-invariant piecewise Euclidean metric of non-positive curvature whenever $P$ is a flag polytope (see~\cite{Gro87}). However, it is not true that $\RZ^n_P$ is always a URC-manifold if $P$ is a flag polytope. Indeed, the direct product of flag polytopes is again a flag polytope, and $\RZ_{P_1\times P_2}^{n_1+n_2} =\RZ_{P_1}^{n_1}\times\RZ_{P_2}^{n_2}$, but a result of Kotschick and L\"oh~\cite{KoLo08} implies that the product of two manifolds of positive dimensions cannot be a URC-manifold. 

\begin{conjecture}
Let $P$ be a flag simple polytope which is not combinatorially equivalent to a product of two simple polytopes of positive dimensions. Then $\RZ_P^n$ is a URC-manifold.
\end{conjecture}

Proposition~\ref{propos_obvious} implies that Conjecture~\ref{conj_main} can be formulated as follows.

\begin{conjecture}
In each dimension $n\ge 2$, there is a  hyperbolic URC-manifold. 
\end{conjecture}

This conjecture remains open for~$n\ge 5$.

For a URC-manifold~$M^n$, let $c_1(M^n)$ and $c_2(M^n)$ be the greatest and the smallest constants respectively such that inequality~\eqref{eq_ineq} holds for every~$X$ and~$z$. It is easy to show that $c_1(M^n)\le \|M^n\|^{-1}\le c_2(M^n)$. It may appear to be interesting to study $c_1$ and~$c_2$ as invariants of URC-manifolds.

The author is grateful to V. M. Buchstaber for fruitful discussions and comments.

\end{document}